\documentclass[aop,MSNbibl,seceqn,dvips]{arximspdf}
\usepackage{accents}

\doi{10.1214/12-AOP812} 
\volume{41}
\issue{5}
\pubyear{2013}
\firstpage{3140}
\lastpage{3156}

\makeatletter
\newcommand{\rrvert}{\vert}
\newcommand{\llvert}{\vert}

\newcommand{\ushort}[1]{\underaccent{\bar}{#1}}

\newcommand{\cal}{\mathcal}

\newtheorem{Lemma}{Lemma}[section]
\newtheorem{Theorem}[Lemma]{Theorem}

\newproclaim{Definition}{Definition}[section]

\makeatother

\begin{document}
\begin{frontmatter}

\title{The empirical cost of optimal incomplete transportation}
\runtitle{Empirical incomplete transportation}

\begin{aug}
\author[A]{\fnms{Eustasio} \snm{del Barrio}\corref{}\thanksref{t1}\ead[label=e1]{tasio@eio.uva.es}}
\and
\author[A]{\fnms{Carlos} \snm{Matr\'an}\thanksref{t1}\ead[label=e2]{matran@eio.uva.es}}
\runauthor{E. del Barrio and C. Matr\'an}
\affiliation{IMUVA, Instituto de Matem\'aticas, Universidad de Valladolid}
\address[A]{Facultad de Ciencias\\
Universidad de Valladolid \\
Paseo de Bel\'en 9 \\
47011 Valladolid\\
Spain \\
\printead{e1}\\
\hphantom{E-mail: }\printead*{e2}} 
\end{aug}

\thankstext{t1}{Supported in part by the Spanish Ministerio de Educaci\'
on y Ciencia and FEDER, Grant
MTM2011-28657-C02-01.}

\received{\smonth{7} \syear{2011}}
\revised{\smonth{9} \syear{2012}}

%
\begin{abstract}
We consider the problem of optimal incomplete transportation between
the empi\-rical measure on an i.i.d. uniform sample on the
$d$-dimensional unit cube $[0,1]^d$ and the true measure. This is a
family of problems lying in between classical optimal transportation
and nearest neighbor problems. We show that the empirical cost of
optimal incomplete transportation vanishes at rate $O_P(n^{-1/d})$,
where $n$ denotes the sample size. In dimension $d\geq3$ the rate is
the same as in classical optimal transportation, but in low dimension
it is (much) higher than the classical rate.
\end{abstract}

%

\end{frontmatter}

\section{Introduction}\label{intro}
Consider two probability measures on $\mathbb{R}^d$, $P$ and $Q$, and
the set $\mathbb{T} (P,Q)$ of maps transporting $P$ to $Q$,
that is, the set of all measurable maps $T\dvtx  \mathbb{R}^d \to\mathbb
{R}^d$ such that, if the initial space is
endowed with the probability $P$, then the distribution of the random
variable $T$ is $Q$.
Monge's optimal transportation problem consists of relocating a certain
amount of mass from its
original distribution to a different target distribution minimizing the
transportation cost.
In more abstract terms, the problem
consists of finding a transportation map $T_0 \in\mathbb{T} (P,Q)$
such that
\[
T_0:=\mathop{\arg\min}_{T \in\mathbb{T} (P,Q)}\int_{\mathbb{R}^d} \bigl\|x
- T(x)\bigr\|^p P(dx).
\]
Here, and throughout the paper, we assume $p\geq1$. Remarkably, under
some smoothness
assumptions, Monge's problem is intimately related to the
$L_p$-Wasserstein distance by
\[
\mathcal{W}_p(P,Q)=\min_{T \in\mathbb{T} (P,Q)} \biggl( \int
_{\mathbb
{R}^d} \bigl\|x - T(x)\bigr\|^p P(dx) \biggr)^{1/p},
\]
where the $L_p$-Wasserstein distance between $P$, $Q$ is defined as
\[
\mathcal{W}_p(P,Q):= \biggl(\inf_{\tau\in\mathcal{M}(P,Q)} \biggl\{
\int\| x-y\|^p \,d\tau(x,y) \biggr\} \biggr)^{1/p},
\]
and $\mathcal{M}(P,Q)$ is the set of probability measures on $\mathbb
{R}^d \times\mathbb{R}^d $ with marginals $P$ and $Q$.
This functional is related to very important problems in mathematics,
the study of which has led to deep
developments in several fields of research and applications, linked to
such important
names as Amp\`ere, Kantorovich, Rubinstein, Zolotarev and Dobrushin,
among others.
To avoid a huge number of references, we refer to the books by
Rachev and R\"uschen\-dorf~\cite{Rachev} and by Villani~\cite{Villani,Villani2},
for an updated account of the interest and implications of the problem.
However, we emphasize the importance of the topic in the development of the
theory of probability metrics and its implications in statistics,
particularly in
goodness of fit problems. Focusing on such kind of problems, the
functional of
interest is $\mathcal{W}_p(P_n,Q)$, where $P_n$ is the empirical
measure associated
to a sample $X_1,\ldots,X_n$, of independent identically distributed
(i.i.d.) random
vectors with law $P$, and $Q$ is any target probability measure on
$\mathbb{R}^d$, or
$\mathcal{W}_p(P_n,Q_n)$, where $Q_n$ stands for the empirical measure
on a second, independent i.i.d.
random sample, $Y_1,\ldots,Y_n$. These empirical versions are connected
to a combinatorial optimization
problem, namely the \textit{optimal matching} problem. In fact,
$\mathcal{W}_p^p(P_n,Q_n)=T_p(n)$, where
%
\begin{equation}
\label{transportfunctional} T_p(n):= \min_\pi
\frac1 n \sum_{i=1}^n\|X_i-Y_{\pi(i)}
\|^p
\end{equation}
and $\pi$ ranges over the permutations of the set $\{1,\ldots,n\}$.
A lot of work has been devoted to analyzing the rate and mode of
convergence of (\ref{transportfunctional}) and several variants of it,
beginning with the seminal paper by Ajtai, Komlos and Tusnady~\cite{AKT}
in the case in which both samples come from the same underlying
probability law $P$. The problem can be equivalently formulated in
terms of
$\mathcal{W}_p(P_n,P)$, the distance between the empirical and true
distributions.
Further references will be provided later, but now let us mention the
series of papers authored by Talagrand~\cite{Talagrand92,Talagrand94}, Talagrand
and Yukich~\cite{Talagrand93} and Dobri\'c and Yukich~\cite{Dobric},
which in
the case when $P$ is the uniform distribution on the $d$-dimensional
unit cube, $[0,1]^d$,
essentially shows that
%
\begin{equation}
\label{untrimmed} \bigl(T_p(n)\bigr)^{1/p}= \cases{
O_P \bigl(n^{-1/2} \bigr), &\quad if $d=1$,
\vspace*{2pt}\cr
\displaystyle O_P
\biggl( \biggl(\frac{\log n} n \biggr)^{1/2} \biggr), &\quad if $d=2$,
\vspace*{2pt}\cr
O_P \bigl(n^{-1/d} \bigr), &\quad if $d\geq3$.}
\end{equation}

This paper deals with the empirical cost of optimal \textit{incomplete}
or \textit{partial} transportation. This is the case in which
the amount of mass required in the target distribution is smaller than
that in the original one.
Then, we do not have to move all the original mass, but we can dismiss
a fraction of it. Of course, we would like to complete
this task with a minimal cost. A more general version is possible if we
admit that we only have to fulfill a fraction
of the target distribution. The general formulation of this problem,
with quadratic cost, has been introduced by
Caffarelli and McCann~\cite{CaffarelliMcCann}, relating it to a
Monge--Amp\`{e}re double obstacle problem.
They obtain remarkable results on the existence, uniqueness and
regularity of the optimal solutions
in a well-separated situation. Figalli~\cite{Figalli} improved the
results covering the case
of nondisjoint supports for the involved probability measures.
Independently, \'Alvarez-Esteban et al.~\cite{Pedro2} introduced the
problem in the context of similarity of probabilities,
obtaining a more general result of existence and uniqueness of the
optimal solution. Moreover~\cite{Pedro2} includes
almost sure consistency of sample solutions to the true ones. In a
subsequent paper, \'Alvarez-Esteban et al.~\cite{Pedro3}
noticed the faster rate of decay of the cost of empirical incomplete
transportation (in the $L_2$ case) and
introduced a procedure for testing similarity of probabilities based on
this fact.

A convenient mathematical formulation of this optimal incomplete
transportation problem can be done with the
help of the concept of \textit{trimmings} of a probability.
%
\begin{Definition} \label{definition1}
Given $0\leq\alpha\leq1$ and Borel probability measures $P$, $R$
on~$\mathbb{R}^d$, we say that $R$ is an $\alpha$-trimming of $P$ if
$R$ is absolutely continuous with respect to $P$, and the
Radon--Nikodym derivative satisfies $\frac{dR}{dP}\leq
\frac{1}{1-\alpha}$. The set of all $\alpha$-trimmings of $P$ will be
denoted by $\mathcal R_\alpha(P)$.
\end{Definition}
Note that in the extreme case $\alpha=0$, $\mathcal R_0 (P)$ is just
$P$, while
$\mathcal R_1 (P)$ is the set of all probability measures absolutely continuous
with respect to $P$. See~\cite{Pedro2} for useful alternative characterizations
of trimmings of a probability, as well as mathematical properties of
the set
$\mathcal R_\alpha(P)$. Turning back to the partial mass
transportation problem,
we could represent the target distribution by the probability $Q$ and
the initial
distribution of mass by $\frac{1}{1-\alpha}P$, $P$ being another
probability if the
mass required in the target distribution is $1-\alpha$ times the mass
in the original locations.
An incomplete transportation plan is then a probability measure $\tau$
on $\mathbb{R}^d\times\mathbb{R}^d$
with second marginal equal to $Q$ and first marginal in $\mathcal
{R}_\alpha(P)$, and the cost
of optimal incomplete transportation is
\[
\mathcal{W}_p\bigl(\mathcal R_\alpha(P),Q\bigr):=\min
_{R\in\mathcal{R}_\alpha
(P) }\mathcal{W}_p(R,Q).
\]
In the more general case, with slackness in the target distribution,
the optimal incomplete transportation cost
would be $\mathcal{W}_p(\mathcal R_{\alpha_1} (P),\mathcal R_{\alpha_2}
(Q))$, the minimal $\mathcal{W}_p$
distance between trimmings of $P$ and $Q$.

This paper gives exact rates of convergence for empirical versions of the
optimal incomplete transportation cost. As with classical optimal
transportation,
the results can be considered in terms of a combinatorial optimization
problem, that we call
\textit{optimal incomplete matching}. To be precise, assume that we can
trim (eliminate) a
fixed proportion $\alpha$ of $X$'s points and also of $Y$'s points, and
we should only search
for the best matching between the nontrimmed samples. Taking for
simplicity $m:=n-\alpha n$ to be
an integer, the new functional of interest is
%
\begin{equation}
\label{newfunctional} T_{p,\alpha}(n): = \min_{X^*,Y^*} \min
_\pi\frac1 m \sum_{j=1}^m
\bigl\| X_j^* -Y^*_{\pi(j)} \bigr\|^p,
\end{equation}
where $\pi$ varies in the set of permutations of $\{1,\ldots,m\}$, $\{
X_1^*,\ldots,X^*_m\}$
ranges in the subsets of size $m$ of $\{X_1,\ldots,X_n\}$, and similarly
$\{Y_1^*,\ldots,Y^*_m\}$ ranges in the subsets of size $m$ of $\{
Y_1,\ldots,Y_n\}$. It is easy to check
that
%
\begin{equation}
\label{partialmatching} T_{p,\alpha}(n)=\mathcal{W}_p^p
\bigl(\mathcal R_\alpha(P_n),\mathcal R_\alpha(Q_n)
\bigr).
\end{equation}
In fact, $\mathcal{W}_p^p(\mathcal R_\alpha(P_n),\mathcal R_\alpha
(Q_n))$ equals the minimum [in $(\pi_{i,j})$]
of the linear function $\sum_{i,j=1}^n \| X_i -Y_{j} \|^p \pi_{i,j}$
subject to the linear constraints
$\sum_{i=1}^n \pi_{i,j}\leq\frac1 {n(1-\alpha)}=\frac1 m$, $\sum
_{j=1}^n \pi_{i,j}\leq\frac1 m$,
$\sum_{i,j=1}^n \pi_{i,j}=1$, $\pi_{i,j}\geq0$ [we are assuming
$m=n(1-\alpha) $]. Rescaling we see that
$m\mathcal{W}_p^p(\mathcal R_\alpha(P_n),\mathcal R_\alpha(Q_n))$ equals
the minimum [in $(x_{i,j}), (a_i),(b_j)$]
of the linear function $\sum_{i,j=1}^n \| X_i -Y_{j} \|^p x_{i,j}$
subject to the linear constraints
$\sum_{i=1}^n x_{i,j}=b_j$, $\sum_{j=1}^n \pi_{i,j}=a_i$, $0\leq
a_i\leq1$,
$0\leq b_j \leq1$, $\sum_{i=1}^n a_i=\sum_{j=1}^n b_j=m$, $x_{i,j}\geq
0$. The constraint matrix in this last
linear program is totally unimodular (see, e.g., Theorem 13.3 in \cite
{PapadimitrouSteiglitz}) and the right-hand
side is integer. Hence, the minimum is attained at some integer
solution, that is,
satisfying $a_i,b_j\in\{0,1\}$, and this implies (\ref{partialmatching}).

We will show the somewhat unexpected result that, independently of the
value of $\alpha\in(0,1)$,
the rates in (\ref{untrimmed}) change to
%
\begin{equation}
\label{trimmed} \bigl({T_{p,\alpha}(n)} \bigr)^{1/p}=O_P
\bigl(n^{-1/d} \bigr)
\end{equation}
for \textit{any} dimension $d\geq1$.
In fact, (\ref{trimmed}) follows from the triangle inequality and
%
\begin{equation}
\label{theproblem} \mathcal{W}_p\bigl(\mathcal R_\alpha(P_n),P
\bigr)=O_P \bigl(n^{-1/d} \bigr),
\end{equation}
which is the formulation we choose for the results we prove. Our
approach relies only on elementary or rather classical tools. In
particular, we do not use subadditivity arguments as in~\cite{Dobric}.
Subadditivity yields a.s. convergence to a constant, rather than just a
rate of convergence. On the other hand, the approach in~\cite{Dobric}
relies on showing subadditivity of a certain Poissonization of the
matching functional (subadditivity does not hold for the original
matching functional; see Remark~1.1 in~\cite{Dobric}). We could also use
that approach here, at least for $p=1$ (otherwise duality for optimal
matching, which is essential in the cited approach, becomes harder to
deal with) but the approximation rate for the Poissonization of the
incomplete matching functional would not allow us to recover the
present result in dimension $d\leq2$.

The study of the rate of convergence of $D(\mathcal R_\alpha(P_n),P)$
for $\alpha\in(0,1)$,
for some probability metrics $D$ was started in del Barrio and Matr\'an
\cite{Tasio}.
In the case of the Wasserstein metric and dimension $d=1$, the results
in~\cite{Tasio} already show
a very different behavior with respect to the untrimmed case, namely
%
\begin{equation}
\label{rate}\quad \frac{n}{(\log n)^\nu}\mathcal{W}_p \bigl(
\mathcal{R}_\alpha(P_n), P\bigr)\to_P 0 \qquad\mbox{for any } \nu>1
\mbox{ and every } \alpha>0.
\end{equation}

We close this \hyperref[intro]{Introduction} mentioning the connection of empirical
optimal incomplete tranportation
to another important problem in probability, that of \textit{random
quantization}.
Taking $\alpha=1$ (full trimming) we have that $\mathcal R_1 (P_n)$ is
the set of all probability measures
concentrated on the sample points, and $\mathcal{W}_p(\mathcal R_1
(P_n),P)$ is the minimal $L_p$-cost
of relocating a mass distributed according to some probability measure
$P$ to a collection
of randomly chosen spots $X_1,\ldots,X_n$. When $X_1,\ldots,X_n$ are
$\mathbb{R}^d$-valued i.i.d. random vectors and $P$ is absolutely continuous,
the problem can be formulated, in Monge's way, as the minimization of
%
\begin{equation}
\label{equacion1} \int_{\mathbb{R}^d}\bigl\|x-\varphi(x) \bigr\|^p
\,dP(x),
\end{equation}
where $\varphi$ varies in the set of all measurable functions with
values in $\{ X_1,\ldots,\break X_n\}$.
Since, for a fixed $x$ in the integrand in (\ref{equacion1}), the
distance $\|x-\varphi(x) \|$ is minimized for
$\varphi(x)=\arg\min_i \|x-X_i\|$, without any constraint on the
capacity to be stored at $X_i$, we obtain that the optimal
$\varphi$ is given by this last expression, and hence the optimal cost equals
%
\begin{equation}
\label{equacion2} \mathcal{W}_p^p\bigl(\mathcal
R_1 (P_n),P\bigr)=\int_{\mathbb{R}^d}\min
_{1\leq i
\leq n}\|x-X_i \|^p \,dP(x).
\end{equation}
Random quantization is a well-studied problem; see, for example, the
Graf and Luschgy monograph
\cite{Graf} or the more recent paper by Yukich~\cite{Yukich}. In particular,
the asymptotic behavior of the $L_p$ quantization error is known, hence
the rate
at which $\mathcal{W}_p(P_{n,1},P)$ vanishes. A trivial consequence of
Definition~\ref{definition1}
is that $\mathcal{R}_{\alpha_1} (P)\subset\mathcal{R}_{\alpha_2} (P)$
if $0\leq\alpha_1\leq\alpha_2\leq1$. This implies
%
\begin{equation}
\label{cotaprincipal} \mathcal{W}_p(\mathcal P_n,P)=
\mathcal{W}_p\bigl(\mathcal R_0 (P_n),P\bigr)
\geq\mathcal{W}_p\bigl(\mathcal R_\alpha(P_n),P
\bigr)\geq\mathcal{W}_p\bigl(\mathcal R_1
(P_n),P\bigr).\hspace*{-28pt}
\end{equation}
Hence rates of convergence for the random quantization error are a
lower bound for rates of convergence of $\mathcal{W}_p(\mathcal
R_\alpha(P_n),P)$
for general $\alpha$. In a first look, classical optimal transportation
is a \textit{global} problem while random quantization is a \textit
{local} one: a point $x$ is mapped through
the optimal map $\varphi$ to a sample point which in the case of random
quantization, is determined just by sample points which are close to $x$
(the nearest neighbor in fact) while in the case of optimal
transportation, two samples with the same sample points in a
neighborhood of $x$ may
result, however, in very different destinations for $x$ due to capacity
constraints. It turns out though, that this different character is only
apparent,
in terms of rates, in dimensions $d=1$ or~$2$. The most relevant fact
which we show in this paper is
that, again in terms of rates, optimal incomplete transportation shows
the same \textit{local} nature as
random quantization in \textit{any} dimension.

The remainder of this paper is organized as follows.
In Section~\ref{preliminary} we give a quick survey on known results
about rates of convergence for optimal
transportation and random quantization. Section~\ref{dimension1}
contains new results for optimal incomplete transportation.
We consider first the case $d=1$, and in this case we construct upper
and lower envelopes for the optimal solution
to the incomplete transportation problem. These are not optimal, but
attain the correct rate.
Finally, we construct a nearly optimal solution in general dimension
starting from the one-dimensional construction.

We will use $EX$ to denote the expected value of a random variable $X$.
By $P(\cdot|B)$ [resp., $E(\cdot| B)$, or even $E_B$] we refer to the
conditional probability (resp., conditional expectation) given the set
$B$. The indicator function of $B$ will be denoted by~$I_B$, while the
notation $\delta_x$ will be reserved for Dirac's probability measure on
the point $x$.
Unless otherwise stated, the random vectors will be assumed to be
defined on the same probability space $(\Omega, \sigma, \nu)$. We write
$\ell^d$ for Lebesgue measure on the space $({\mathbb{R}^d}, \beta)$.
Finally, convergence in probability (resp., weak convergence of
probabilities) will be denoted by $\to_p$ (resp., by $ \to_{w} $), and
${\cal L}(X)$ will denote the law of the random vector $X$.

\section{Preliminary results}\label{preliminary}

The results on the asymptotic behavior of $L_p$-Wasserstein distances
between the empirical
and parent distributions in the one-dimensional case have been obtained
through a quantile
representation. If $F$ and $G$ are the distribution functions of $P$
and $Q$ and $F^{-1}$ and
$G^{-1}$ are the
respective quantile functions, then (see, e.g., Bickel and Freedman
\cite{Bickel})
%
\begin{equation}
\label{prop1} {\cal W}_p(P,Q)= \biggl[\int_0^1
\bigl(F^{-1}(t)-G^{-1}(t) \bigr)^p \,dt
\biggr]^{1/p}
\end{equation}
[where
$F^{-1}(t)=\inf
\{s\dvtx  F(s) \geq t\}$].
In particular, when $P$ is the uniform distribution on $(0,1)$, this
representation leads to
%
\begin{equation}
\label{BGU} \sqrt{n} \mathcal{W}_p(P_{n},P)
\to_w \biggl[ \int_0^1
\bigl({B(t)} \bigr)^p\,dt \biggr]^{1/p}
\end{equation}
with $B(t)$ a Brownian Bridge on $[0,1]$; see, for example, \cite
{CsorgoHorvath}.

For dimension $d>1$, there are not explicit expressions for the optimal
transportation maps, and limit distribution results as in (\ref{BGU})
are not available. Rates of convergence to 0 of $\mathcal{W}_p(P_n,P)$
can be given based on different approaches. The case $d=2$ is the most
interesting from the point of view of the mass transportation problem.
Ajtai, Koml\'os and Tusn\'adi~\cite{AKT} showed that, with probability
$1-o(1)$,
\[
C_1 \biggl(\frac{\log n}n \biggr)^{1/2} < \mathcal
{W}_1(P_n,Q_n)< C_2 \biggl(
\frac{\log n}n \biggr)^{1/2},
\]
where $P_n$ and $Q_n$ are the sample distributions corresponding to two
independent
samples obtained from the uniform distribution on the unit square, $U([0,1]^2)$.
Their combinatorial partition scheme method was refined in Talagrand
and Yukich~\cite{Talagrand93}
to show [Theorem 1 and Remark (ii) there] that for some constant, $C(p)$,
%
\begin{equation}
\label{AKTbound} E\bigl(\mathcal{W}_p\bigl(P_n,U
\bigl([0,1]^2\bigr)\bigr)\bigr)\leq C(p) \biggl( {\frac{\log
n}n}
\biggr)^{1/2}.
\end{equation}
The case $d\geq3$ (and uniform distribution on the $d$-dimensional
unit cube) is covered in
Talagrand~\cite{Talagrand94}. That paper uses a different approach,
based on duality for the optimal
transportation problem to give a result (Theorem~1.1), formulated for a
very general class of costs
functions which includes exponential costs and, as a consequence, implies
%
\begin{equation}
\label{Talagrandbound} E\bigl(\mathcal{W}_p\bigl(P_n,U
\bigl([0,1]^d\bigr)\bigr)\bigr)\leq C(k,p) \frac1 {n^{1/d}}.
\end{equation}
Further results, dealing with distributions other than the uniform,
possibly with unbounded support, are
given in Barthe and Bordenave~\cite{BartheBordenave}.

As we already noted in the \hyperref[intro]{Introduction}, the so-called random
quantizers provide an easy way of
giving a lower bound for the rates of convergence of our interest. We
give a simple version of
the mean asymptotics for the random quantizers, rewritten in terms of
the $L_p$-Wasserstein distance
between the set of ``fully trimmed'' sample probabilities and the
theoretical distribution, that
suffices for our purposes; this is a particular case of Theorem 9.1 in
\cite{Graf}.
%
\begin{Theorem}\label{tasaconrecorte1} If $X_1,\ldots,X_n$ are i.i.d.
random vectors uniformly
distributed on $[0,1]^d$, then
\[
n^{p/d}E\bigl(\mathcal{W}_p^p\bigl(R_1
(P_n),U\bigl([0,1]^d\bigr)\bigr)\bigr)\to\Gamma{ \biggl(1+
\frac p k \biggr){\omega_d^{-p/d}}} \qquad\mbox{as } n\to\infty,
\]
where $\omega_d=\frac{\pi^{d/2}} {\Gamma(1+ d/2 )}$.
\end{Theorem}


For optimal incomplete transportation, a first result on rates of
convergence is Theorem 5 in the Appendix of~\cite{Pedro3}, for
dimension 1, but it has been largely improved in~\cite{Tasio}, in the
terms expressed in (\ref{rate}). For dimension 2 our approach was not
successful in going beyond the characteristic ``$\log n$'' term in the
Ajtai--K\'omlos--Tusn\'ady result (\ref{untrimmed}). This task and the
fundamental improvement in dimension 1 are the main goals in this
paper. Moreover we notice del Barrio and Matr\'an~\cite{Tasio} also
treat the improvement of the ``in probability bounds'' involved in
(\ref{theproblem}) to almost surely bounds. This follows Talagrand's
approach~\cite{Talagrand92}, continued by Dobri\'c and Yukich \cite
{Dobric}, but, in our case, using a powerful concentration inequality
of Boucheron et~al.~\cite{Boucheron}.\looseness=1

\section{Rates of convergence}\label{dimension1}

We focus first on the one-dimensional case. Let us consider $n$
distinct points $x_1< \cdots<x_n\in(0,1)$
and set $P_n=\frac1 n \sum_{i=1}^n \delta_{x_i}$, where $\delta_x$
denotes Dirac's measure on $x$. An $\alpha$-trimming of $P_n$ can be written,
in terms of a vector $h=(h_1,\ldots,h_{n-1})$, as $(P_n)_h=\sum_{i=1}^n
b_i \delta_{x_i}$ with $0\leq b_i=h_i-h_{i-1}\leq\frac1{n(1-\alpha)}$
(we set, for convenience, $h_0=0, h_n=1$). We therefore write
%
\begin{eqnarray}
\label{setCa}
&&\mathcal{C}_{\alpha, n}:= \biggl\{h=(h_1,\ldots,h_{n-1})\in\mathbb{R}^{n-1}\dvtx  0\leq h_i-h_{i-1}
\leq\frac1{n(1-\alpha)},\nonumber\\[-8pt]\\[-8pt]
&&\hspace*{212pt}\qquad i=1,\ldots, n \biggr\}.\nonumber
\end{eqnarray}
Our first result is an elementary, but useful, representation of
$\mathcal{W}_p(\mathcal{R}_\alpha(P_n),P)$.

\begin{Lemma}\label{expresionutil}
If $x_1< \cdots<x_n\in(0,1)$, $P_n=\frac1 n \sum_{i=1}^n \delta
_{x_i}$, $P$ is the uniform distribution on $[0,1]$, $\mathcal
{C}_{\alpha, n}$
is defined by (\ref{setCa}) and $p\geq1$, then
\begin{eqnarray*}
&&
\mathcal{W}_p^p\bigl(\mathcal{R}_\alpha(P_n),P
\bigr)\\
&&\qquad=\frac1 {p+1} \bigl(x_1^{p+1}+(1-x_n)^{p+1}
\bigr)
\\
&&\qquad\quad{}+\frac1 {2^p (p+1)}\sum_{i=1}^{n-1}(x_{i+1}-x_i)^{p+1}
\\
&&\qquad\quad{}+ \min_{h\in\mathcal{C}_{\alpha,n}} \frac1 {p+1} \sum_{i=1}^{n-1}
\biggl(\frac{x_{i+1}-x_i}{2} \biggr)^{p+1} f_p \biggl(
\frac{h_i-({x_{i+1}+x_i})/{2}}{({x_{i+1}-x_{i}})/{2}} \biggr),
\end{eqnarray*}
where $f_p(y)=(1+|y|)^{p+1}+ (1-|y|)^{(p+1)}-2$ and $t^{(p)}$ denotes
the odd extension to $(-\infty,\infty)$ of
the function $t^p$ on $[0,\infty)$.
\end{Lemma}

\begin{pf}
We note first that the quantile function
associated to $(P_n)_h$
takes the value $x_i$ in the interval $(h_{i-1},h_i]$. Hence, using
(\ref{prop1}),
we see that
$\mathcal{W}_p^p(\mathcal{R}_\alpha(P_n), P)=\min_{h\in\mathcal
{C}_{\alpha,n}} \sum_{i=1}^n A_i$,
where $A_i=\int_{h_{i-1}}^{h_i} |x_i-t|^p \,dt$. Since\vadjust{\goodbreak} $t^{(p+1)}/(p+1)$
is a primitive of
$|t|^p$, we can write
\begin{eqnarray*}
A_i &=&\int_{h_{i-1}}^{x_i}
|x_i-t|^p \,dt+\int_{x_i}^{h_i}
|x_i-t|^p \,dt
\\
&=&\frac1{p+1} \bigl[(x_i-h_{i-1})^{(p+1)}+
(h_i-x_i)^{(p+1)} \bigr].
\end{eqnarray*}
From this, recalling that $h_0=0$, $h_1=1$, we get
\[
\sum_{i=1}^n A_i= \frac1
{p+1} \Biggl[ \bigl(x_1^{p+1}+(1-x_n)^{p+1}
\bigr)+\frac1 {2^p }\sum_{i=1}^{n-1}(x_{i+1}-x_i)^{p+1}+
\sum_{i=1}^{n-1} B_i \Biggr]
\]
with
$B_i= (x_{i+1}-h_i)^{(p+1)}+(h_i-x_i)^{(p+1)}-2 (\frac
{x_{i+1}-x_i}{2} )^{p+1}$.
Now it is easy to see that $B_i= (\frac{x_{i+1}-x_i}{2}
)^{p+1} f_p
(\frac{h_i-({x_{i+1}+x_i})/{2}}{({x_{i+1}-x_{i}})/{2}}
)$, which completes the proof.
\end{pf}

The function $f_p$ in Lemma~\ref{expresionutil} is a piecewise
polynomial for integer $p$. For instance
$f_1(y)=2y^2$, $|y|\leq1$, $f_1(y)=2(2|y|-1)$, $|y|> 1$;
$f_2(y)=6y^2$, $y\in\mathbb{R}$. For general $p\geq1$,
$f_p$ is a nonnegative, even and convex function, strictly increasing
on $[0,\infty)$, which attains its minimum
at $y=0$, with $f_p(0)=0$. This suggests that a good trimming vector
$h=(h_1,\ldots,h_{n-1})\in\mathcal{C}_{\alpha,n}$
should be as close as possible to the midranks, ${\frac
{x_{i}+x_{i+1}}{2}}$. With this observation in mind, we denote
\[
\hat{h}=\mathop{\arg\min}_{h\in\mathcal{C}_{\alpha, n}} \sum_{i=1}^{n-1}
\biggl( \frac{x_{i+1}-x_i}{2} \biggr)^{p+1} f_p \biggl(
\frac{h_i-({x_{i+1}+x_i})/{2}}{
({x_{i+1}-x_{i}})/{2}} \biggr)
\]
and define
\[
u_i=\max_{i\leq j\leq n-1} \biggl({ \frac{x_j+x_{j+1}}{2} -
\frac{1}{1-\alpha} \frac j n } \biggr)\vee{ \frac{-\alpha
}{1-\alpha}},\qquad i=1,\ldots,n-1,
\]
$u_n=-\frac\alpha{1-\alpha}$, $\bar{f}_0=0$, and $\bar{f}_i=u_i
\wedge0, i=1,\ldots, n$.
Finally, we set
\[
\bar{h}_i=\bar{f}_i+\frac1 {1-\alpha} \frac i n,\qquad i=0,\ldots,n.
\]
Note that, for any $h=(h_1,\ldots,h_{n-1})\in\mathcal{C}_{\alpha,n}$,
$h_i-\frac i {n(1-\alpha)}$ is
a sequence that decreases from 0 to $-\frac\alpha{1-\alpha}$,
while $\bar{f}_i$ is the lowest decreasing sequence from 0 to $-\frac
\alpha{1-\alpha}$
which lies above the sequence $\frac{x_i+x_{i+1}}{2} -\frac
{i}{n(1-\alpha)}$.
In the next result we see that $\bar{h}_i$ is a feasible trimming and
that feasible solutions that exceed this one cannot be optimal.
%
\begin{Lemma}\label{splitting}
$\bar{h}_0=0$, $\bar{h}_n=1$ and if $\bar{h}=(\bar{h}_{1},\ldots,\bar
{h}_{n-1})$,
then $\bar{h}\in\mathcal{C}_{\alpha,n}$. Furthermore,
\[
\hat{h}_i\leq\bar{h}_i,\qquad i=1,\ldots,n-1.
\]
\end{Lemma}

\begin{pf}
$\bar{h}_0=0$ and $\bar{h}_n=1$ are obvious.
To prove $\bar{h}\in\mathcal{C}_{\alpha,n}$ we check, equivalently,
that $0\geq\bar{f}_i-\bar{f}_{i-1}\geq-\frac1{1-\alpha}\frac1 n$,
$i=1,\ldots,n$. Clearly, $u_1\geq\cdots\geq u_n$ and $\bar{f}_1\leq
0=\bar{f}_0$ and, consequently,
$\bar{f}_0\geq\cdots\geq\bar{f}_{n}$. To see that $\bar{f}_i-\bar
{f}_{i-1}\geq-\frac1{n(1-\alpha) }$ observe that $u_i=u_{i-1}$
unless $u_{i-1}=\frac{x_{i-1}+x_i}{2}-\frac{i-1}{n(1-\alpha)}$, but then
\begin{eqnarray*}
u_i-u_{i-1}&\geq& \biggl(\frac{x_{i}+x_{i+1}}{2}-
\frac
{i}{n(1-\alpha)} \biggr) - \biggl(\frac{x_{i-1}+x_i}{2}-\frac
{i-1}{n(1-\alpha)} \biggr)
\\
&=&\frac{x_{i+1}-x_{i-1}}{2}-\frac1{n(1-\alpha)}\\
&\geq&-\frac1{n(1-\alpha)}
\end{eqnarray*}
and the claim follows. We show now that $\hat{f}_i\leq\bar{f}_i$,
where $\hat{f}_i=\hat{h}_i-\frac1 {1-\alpha} \frac i n$.
Since $u_i\geq\frac{x_i+x_{i+1}}{2}-\frac{1}{1-\alpha} \frac i n$, we
see that $\bar{f}_i \geq\frac{x_i+x_{i+1}}{2}-\frac{1}{1-\alpha} \frac
i n$,
provided $\frac{x_i+x_{i+1}}{2}-\frac{1}{1-\alpha} \frac i n\leq0$.
Now, if $u_i\geq0$, then $\bar{f}_i=0\geq\hat{f}_i$. Let us assume that
$u_i<0$ (hence $\bar{f}_i \geq\frac{x_i+x_{i+1}}{2}-\frac{1}{1-\alpha}
\frac i n$) and $\bar{f}_j\geq\hat{f}_j$, $j<i$, but $\bar{f}_i< \hat
{f}_i$. Let us write $k$ for the smallest integer $k>i$ such that $\bar
{f}_k\geq\hat{f}_k$ (observe that $k\leq n$ since $\bar{f}_n= \hat
{f}_n=-\frac\alpha{1-\alpha}$).
We define $\tilde{f}_j=\hat{f}_j$ if $j<i$ or $j\geq k$ and $\tilde
{f}_j=\bar{f}_j$ if $i\leq j<k$. Also, write $\tilde{h}_j=\tilde{f}_j-
\frac1 {1-\alpha} \frac j n$. Clearly, $\tilde{h}\in\mathcal
{C}_{\alpha,n}$. But
for integer $j\in[i,k)$ we have $\hat{h}_j>\tilde{h}_j\geq\frac{x_j
+x_{j+1}}{2}$, which implies
$\llvert\hat{h}_j -{ \frac{x_j+x_{j+1}}{2}} \rrvert>\llvert\tilde
{h}_j -{ \frac{x_j+x_{j+1}}{2}} \rrvert$. Consequently,
\begin{eqnarray*}
&&
\sum_{j=1}^{n-1} (x_{j+1}-x_j)^{p+1}f_p
\biggl(\frac{\tilde{h}_j
-{ ({x_j+x_{j+1}})/{2}}}{({x_{j+1}-x_j})/{2}} \biggr)\\
&&\qquad<\sum_{j=1}^{n-1}
(x_{j+1}-x_j)^{p+1}f_p \biggl(
\frac{\hat{h}_j -{
({x_j+x_{j+1}})/{2}}}{({x_{j+1}-x_j})/{2}} \biggr),
\end{eqnarray*}
against optimality of $\hat{h}$. Hence, $\bar{f}_i\geq\hat{f}_i$ for
all $i$, and the upper bound for $\hat{h}_i$ follows.
\end{pf}

A similar lower bound for $\hat{h}$ can be obtained taking
$\ushort{f}_i$ to be the greatest decreasing sequence from 0 to $-\frac
\alpha{1-\alpha}$ which lies below the sequence $\frac{x_i+x_{i+1}}{2}
-\frac{i}{n(1-\alpha)}$ and setting $\ushort{h}_i=\ushort{f}_i+\frac
{i}{n(1-\alpha)}$. We note\vspace*{1pt} also that Lemma~\ref{expresionutil},
combined with Lem\-ma~\ref{splitting}, gives the following useful lower
and upper bounds for the incomplete transportation cost. To be precise,
%
\begin{eqnarray}
\label{upperlower}\qquad V_n(p) &\leq&\mathcal{W}_p^p
\bigl(\mathcal{R}_\alpha(P_n),P\bigr)
\nonumber\\[-8pt]\\[-8pt]
&\leq&V_n(p)+\frac1 {p+1} \sum
_{i=1}^{n-1} \biggl(\frac
{x_{i+1}-x_i}{2}
\biggr)^{p+1} f_p \biggl(\frac{\bar{h}_i-({x_{i+1}+x_i})/{2}}{
({x_{i+1}-x_{i}})/{2}} \biggr),\nonumber
\end{eqnarray}
where $V_n(p)=\frac1 {p+1} (x_1^{p+1}+(1-x_n)^{p+1} )+\frac
1 {2^p (p+1)}\sum_{i=1}^{n-1}(x_{i+1}-x_i)^{p+1}$. We could replace
$\bar{h}_i$ with
$\ushort{h}_i$ or $\tilde{h}_i=(\bar{h}_i+\ushort{h}_i)/2$, but in
terms of rates, the upper bound above cannot be improved, as we will
see later.

Next, we consider the case of a uniform random sample on the unit
interval, namely $X_1,\ldots,X_n$ are i.i.d. $U(0,1)$ r.v.'s,
$(x_1,\ldots,x_n)=(X_{(1)},\ldots,X_{(n)})$ is the order statistic and
$P_n$ the empirical distribution on the sample. We will use the well-known
fact
%
\begin{equation}
\label{equivalence} (X_{(1)},\ldots,X_{(n)})\stackrel{d} {=}
\biggl(\frac{S_1}{S_{n+1}},\ldots,\frac{S_n}{S_{n+1}} \biggr),
\end{equation}
where, $S_i=\xi_1+\cdots+\xi_i$ and $\{ \xi_i\}_{i=1}^\infty$ are
i.i.d. exponentials random variables with unit mean.
The following elementary lemma about the concentration of the $S_i$'s
around their means will be
used repeatedly in the remainder of this section.
%
\begin{Lemma}\label{concentracionSi}
If $t>0$, then
\[
P(S_i-i>t)\leq e^{-t} \biggl(1+\frac t i
\biggr)^i,
\]
while for $0<t<i$
\[
P(i-S_i>t)\leq e^t \biggl(1-\frac t i
\biggr)^i.
\]
\end{Lemma}

\begin{pf}
This is just Chernoff's inequality; see, for
example,~\cite{Massart2007}, page~16.
\end{pf}

We are ready now to give the rate of convergence of $\mathcal
{W}_p(\mathcal{R}_\alpha(P_n),P)$ in the one-dimensional setup.
%
\begin{Theorem} \label{cotaenesperanza} If $P$ is the uniform
distribution on $[0,1]$, $X_1,\ldots,X_n$ are i.i.d. random variables
with common distribution $P$, $P_n$ is the empirical measure on
$X_1,\ldots,X_n$, $\alpha\in(0,1)$ and $p\geq1$, then there exist
constants, $C_p(\alpha)$, depending only on $p$ and $\alpha$, $c_p>0$
depending only on $p$, such that for every $n\geq1$,
\[
\frac{c_p}{n^p}\leq E\bigl(\mathcal{W}^p_p\bigl(
\mathcal{R}_\alpha(P_n),P\bigr)\bigr)\leq\frac{C_p(\alpha)}{n^p}.
\]
\end{Theorem}

\begin{pf}
For the lower bound simply observe that $E( n^p\mathcal{W}^p_p(\mathcal
{R}_\alpha(P_n),P))\geq n^p E(V_n(p))$, with $V_n(p)$ as in
(\ref{upperlower}). The spacing $X_{(i+1)}-X_{(i)}$ follows a beta
distribution with parameters $1$ and $n+1$, and from this fact it follows
that $n^p E(V_n(p))=\frac{n^p \Gamma(n+1)\Gamma(p+2)}{(p+1)\Gamma
(n+p+2)} (2+\frac{n-1}{2^p})$. It is easy to check (using Stirling's
formula, e.g.) that $n^p E(V_n(p))\to\frac{\Gamma
(p+2)}{2^p(p+1)}>0$ as $n\to\infty$, and
hence we can take $c_p= \min_{n\geq1} n^p E(V_n(p))>0$.\vadjust{\goodbreak}

For the upper bound we use the representation (\ref{equivalence}), fix
$\theta\in(1-\alpha,1)$
and write $Z=\mathcal{W}^p_p(\mathcal{R}_\alpha(P_n),P))$. Then
we split $E(Z)$,
\[
E(Z)=E\biggl(ZI\biggl( \frac{S_{n+1}}{n}<\theta\biggr)\biggr)+E\biggl(Z I
\biggl(\frac
{S_{n+1}}{n}\geq\theta\biggr)\biggr):=E(Z_1)+E(Z_2)
\]
and proceed to bound $E(Z_i)$, $i=1,2$. To deal with $E(Z_1)$ we note
that $Z\leq\mathcal{W}^p_p(P_n,P)=\int_0^1 |G_n^{-1}(t)-t|^p \,dt$,
$G_n$ being the distribution function asociated to $P_n$. A simple
computation, similar to the proof of Lemma~\ref{expresionutil}, shows
$\int_0^1 |G_n^{-1}(t)-t|^p\,dt=\int_0^1 |G_n(t)-t|^p \,dt$, both terms
equaling, in fact,
\[
{ \frac1 {p+1}} \sum_{i=1}^n \biggl[\biggl|{
\frac i n -X_{(i)}}\biggr|^{(p+1)}- \biggl|{ \frac{i-1} n
-X_{(i)}}\biggr|^{(p+1)} \biggr].
\]
Hence, from Schwarz's inequality we get
\begin{eqnarray*}
E(Z_1)&\leq& \bigl(E\bigl(Z^2\bigr)\bigr)^{1/2}
P\biggl( \frac{S_{n+1}}{n}<\theta\biggr)^{1/2}
\\
&\leq& \biggl( E \biggl(\biggl( \int^1_0
\bigl|G_n(t)-t\bigr|^p \,dt\biggr) \biggr)^2
\biggr)^{1/2} P\biggl( \frac{S_{n+1}}{n}<\theta\biggr)^{1/2}
\\
&\leq& \biggl(\int^1_0 E\bigl|G_n(t)-t\bigr|^{2p}
\,dt\biggr)^{1/2} P\biggl( \frac{S_{n+1}}{n}<\theta
\biggr)^{1/2}.
\end{eqnarray*}
Using the fact that $P(|G_n(t)-t|>\varepsilon)\leq2 e^{-2n\varepsilon
^2}$ (this follows from Hoeffding's inequality applied to Bernoulli
random variables), we see that
$E(|G_n(t)-t|^{2p})\leq p 2^{1-p} \Gamma(p) n^{-p}$.
Also, from Lemma~\ref{concentracionSi} we get $P( \frac
{S_{n+1}}{n}<\theta)=
P((n+1)-S_{n+1}>n(1-\theta)+1)\leq e^{n(1-\theta)+1} (\frac{n\theta
} {n+1} )^{n+1}$. Combining these two estimates we get
%
\begin{equation}
\label{cotaZ1} E\bigl(n^pZ_1\bigr)\leq\bigl( p
2^{1-p} \Gamma(p) \theta e \bigl(\theta e^{1-\theta}
\bigr)^n n^p \bigr)^{1/2}.
\end{equation}
The last upper bound is a vanishing sequence (hence, bounded) since
$\theta e^{1-\theta}<1$. Note that the bound depends on $\alpha$ through
the choice of $\theta\in(1-\alpha,1)$.

We consider now $E(Z_2)$ and recall (\ref{upperlower}). We have
$ E(n^p V_n(p) I(\frac{S_{n+1}}{n}\geq\theta))\leq
E(n^p V_n(p))\leq\sup_{m\geq1} E(m^p V_m(p))<\infty$ since, as noted
above,\break $E(n^p V_n(p))$ is a convergent sequence as $n\to\infty$.
Observe now that $f_p(y)\leq2^{p+1}-2$ if $|y|\leq1$, while
$f_p(y)\leq2^{p+1}(p+1)|y|^p -2$ if $|y| \geq1$. Therefore
$f_p(y)\leq2^{p+1}(1+(p+1)|y|^p)$, and it suffices to give an upper
bound for
\begin{eqnarray*}
&&
E \Biggl(n^p I\biggl(\frac{S_{n+1}}{n}\geq\theta\biggr)
\sum_{i=1}^{n-1} (X_{(i+1)}-X_{(i)})
\biggl\llvert{ \bar{h}_i} -{ \frac{X_{(i)}+X_{(i+1)}}{2}} \biggr\rrvert
^p \Biggr)
\\
&&\quad\leq\frac1 {\theta^{p+1}} \frac1 n \sum
_{i=1}^{n-1} E\biggl( \xi_{i+1}\biggl\llvert
F_i - \biggl(S_i +\frac{\xi_{i+1}}{2}-\frac1 {1-\alpha}
\frac{S_{n+1}} n i\biggr) \biggr\rrvert^{p} I\biggl(
\frac{S_{n+1}}{n}\geq\theta\biggr)\biggr),
\end{eqnarray*}
where we are using representation (\ref{equivalence}) and
\[
F_i= \biggl( \biggl(\max_{i\leq j\leq n-1}
\biggl(S_j +\frac{\xi
_{j+1}}{2}-\frac1 {1-\alpha} \frac{S_{n+1}} n j
\biggr) \biggr) \vee\biggl(-\frac\alpha{1-\alpha} S_{n+1} \biggr)
\biggr) \wedge0.
\]
It only remains to find an upper bound for $E(U_n)$ with
\[
U_n:=\frac1 n \sum_{i=1}^{n-1}
\xi_{i+1}\biggl\llvert F_i - \biggl(S_i +
\frac{\xi_{i+1}}{2}-\frac1 {1-\alpha} \frac
{S_{n+1}} n i\biggr) \biggr\rrvert
^p I\biggl(\frac{S_{n+1}}{n}\geq\theta\biggr).
\]
We split the sum in $U_n$ into three terms,
$U_n=U_n^{(1)}+U_n^{(2)}+U_n^{(3)}$, $U_n^{(1)}$ collecting the
summands with $F_i=0$, $U_n^{(3)}$ those with
$F_i=-\frac\alpha{1-\alpha} S_{n+1}$ and $U_{n}^{(2)}$ the others. We
bound first $E(U_{n}^{(1)})$.
We write $K=\frac\theta{1-\alpha}$ and note that $K>1$. Now,
\[
U_{n}^{(1)}\leq\frac1 n \sum_{i=1}^{n-1}
\xi_{i+1} \biggl\llvert S_i +\frac{\xi
_{i+1}}{2}-
\frac{1}{1-\alpha} \frac{S_{n+1}} n i \biggr\rrvert^p I\Bigl(\max
_{i\leq j\leq n-1} (S_{j+1} -K j) \geq0\Bigr).
\]
Convexity implies that $E(\frac{S_i} i)^s\leq E\xi_1^s$ for
$s\geq1$. From the Schwarz inequality and the moment inequality
$E|X+Y|^p\leq2^{p-1}(E|X|^p+E|Y|^p)$, $p\geq1$, we get that $(E(\xi
_{i+1} | S_i +\frac{\xi_{i+1}}{2}-\frac1 {1-\alpha} \frac
{S_{n+1}} n i |^p)^2)^{1/2}\leq C_1 i^p$ for
some absolute constant $C_1$ (not depending on $i$ or $n$). On the
other hand, using again Lemma~\ref{concentracionSi} with $i=j+1$ and
$t=Kj-(j+1)$ [which is positive for
$j>(1-\alpha)/(\theta-1+\alpha)$] we have
$P(S_{j+1} -K j \geq0)\leq e^{-(K-1)j+1} (\frac{j
(K-1)}{j+1} )^{j+1}\leq K e q^j$, where \mbox{$q=K e^{-(K-1)}<1$}.\vspace*{1pt} Since
$P(\max_{i\leq j\leq n-1} (S_{j+1} -K j) \geq0)\leq\sum_{j\geq
i}P(S_{j+1} -K j \geq0)$, we get, for some constant, $C_2$,
\[
P\Bigl(\max_{i\leq j\leq n-1} (S_{j+1} -K j) \geq0\Bigr)\leq
C_2 \frac{q^i}{1-q}
\]
[$C_2=Ke$ suffices for $i\geq(1-\alpha)/(\theta-1+\alpha)$; with a
larger constant, if necessary, the bound is true for all $i$].
Combining the last bounds and Schwarz's inequality we obtain, with a
new constant $C_3$,
\[
E\bigl(U_{n}^{(1)}\bigr)\leq\frac{C_3}{\sqrt{1-q}} \frac1 n \sum
_{i=1}^{n-1} i^p q^{i/2}
\]
and, again, the right-hand side is a vanishing sequence.

To deal with $U_{n}^{(3)}$ we define $\xi'_{i}=\xi_{n+2-i}$,
$S_i'=\xi'_1+\cdots+\xi'_i$, $i=1,\ldots,\break n+1$. Observe that
$S_i'=S_{n+1}-S_{n+1-i}$, $i=1,\ldots, n$, and $S_{n+1}'=S_{n+1}$. We
also write $A_i=(\frac{S_{n+1}}{n}\geq\theta,\max_{i\leq j\leq n-1}
(S_j +\frac{\xi_{j+1}}{2}-\frac1 {1-\alpha} \frac{S_{n+1}} n j)\leq
-\frac{\alpha}{1-\alpha} S_{n+1} )$. Then
\begin{eqnarray*}
U_n^{(3)}&=&\frac1 n \sum
_{i=1}^{n-1} \xi_{i+1}\biggl\llvert-
\frac{\alpha}{1-\alpha} S_{n+1} - \biggl(S_i +\frac{\xi_{i+1}}{2}-
\frac1 {1-\alpha} \frac{S_{n+1}} n i\biggr) \biggr\rrvert^p
I_{A_i}
\\
&=&\frac1 n \sum_{j=1}^{n-1}
\xi_{n+1-j}\biggl\llvert-\frac{\alpha
}{1-\alpha} S_{n+1} \\
&&\qquad\hspace*{38pt}{}-
\biggl(S_{n-j} +\frac{\xi_{n+1-j}}{2}-\frac1 {1-\alpha} \frac
{S_{n+1}} n
(n-j)\biggr) \biggr\rrvert^p I_{A_{n-j}}.
\end{eqnarray*}
With the above notation we see that $ -\frac{\alpha}{1-\alpha} S_{n+1}
- (S_{n-j} +\frac{\xi_{n+1-j}}{2}-\break\frac1 {1-\alpha} \frac {S_{n+1}} n
(n-j))=S_j'+\frac{\xi'_{j+1}}2-\frac1 {1-\alpha}\frac {S'_{n+1}}n j$,
while $A_{n-j}=(\frac{S'_{n+1}}{n}\geq\theta,\break\min_{1\leq k\leq j}
(S_k'+\frac{\xi'_{k+1}}2-\frac1 {1-\alpha}\frac{S'_{n+1}}n k)\geq
0)\subset (\frac{S'_{n+1}}{n}\geq\theta,\max_{j\leq k\leq n-1}
(S_k'+\frac{\xi'_{k+1}}2-\frac1 {1-\alpha}\frac{S'_{n+1}}n
k)\geq0):=B_j$. These observations imply that
\begin{eqnarray*}
U_n^{(3)}&=&\frac1 n \sum_{j=1}^{n-1}
\xi'_{j+1}\biggl\llvert S_j'+
\frac{\xi'_{j+1}}2-\frac1 {1-\alpha}\frac{S'_{n+1}}n j \biggr\rrvert
^p I_{A_{n-j}}
\\
&\leq& \frac1 n \sum_{j=1}^{n-1}
\xi'_{j+1}\biggl\llvert S_j'+
\frac{\xi'_{j+1}}2-\frac1 {1-\alpha}\frac{S'_{n+1}}n j \biggr\rrvert
^p I_{B_{j}}.
\end{eqnarray*}
The last upper bound and $U_n^{(1)}$ are equally distributed. Hence
$E(U_{n}^{(3)})\leq
E(U_{n}^{(1)})\to0$.

We turn now to the central part, $U_{n}^{(2)}$. Obviously
%
\begin{equation}
U_{n}^{(2)} \leq\frac1 n \sum_{i=1}^{n-1}
\xi_{i+1} Z_{i}^p,
\end{equation}
where $Z_i= \sup_{j\geq i} ((S_{j+1}-S_{i})- K (j-i))_+$. Once more we
use Schwarz's inequality to get
$E(\xi_{i+1}Z_i^p)\leq(E\xi^2_{i+1})^{1/2}(EZ^{2p}_{i})^{1/2}\le\sqrt
{2} (EZ_0^{2p})^{1/2}$. Thus, it only remains to show $EZ_0^{2p}<\infty
$. Chernoff's
inequality yields
$P(S_i-Ki>t)\leq e^{-(t+(K-1)i)} (K+\frac t i )^i$. From this
and the fact $\int_{0}^\infty e^{-t} t^{l}\,dt=l!$, $l\in\mathbb{N}$, we
get, for integer
$k\geq1$,
%
\begin{eqnarray}\label{cota23}
\int_0^\infty t^k
P(S_i-Ki>t) \,dt &\leq& \int_0^\infty
e^{-t} e^{-(K-1)i} \sum_{j=0}^i
\pmatrix{i
\cr
j}K^{j} \frac{t^{i-j+k}}{i^{i-j}} \,dt
\nonumber\\
&\leq& \frac{e^i i!}{i^i}\sum_{j=0}^i
\frac{e^{-K i}
(Ki)^j}{j!}(i+k-j)^k
\nonumber\\[-8pt]\\[-8pt]
&\leq& (i+k)^k \frac{e^i i!}{i^i} \sum
_{j=0}^i \frac{e^{-K i}
(Ki)^j}{j!}
\nonumber\\
&=& (i+k)^k \frac{e^i i!}{i^i}P(N_{Ki}\leq i),\nonumber
\end{eqnarray}
where $N_\lambda$ denotes a random variable having Poisson distribution
with mean~$\lambda$. Chernoff's inequality (for the left tail) gives
\[
P(N_{Ki}\leq i)=P\bigl(Ki-N_{K i}\geq(K-1)i\bigr)\leq\exp
\biggl(-i K h\biggl( -\frac{K-1}{K}\biggr)\biggr),
\]
where $h(u)=(1+u)\log(1+u) -u$, $u\geq-1$; see, for example, \cite
{Massart2007}, page 19. This, (\ref{cota23}) and the fact $\frac{e^i
i!}{i^i}\leq C \sqrt{i}$
for some constant $C$, imply
\begin{eqnarray*}
E\bigl(Z_0^{k+1}\bigr)&=&(k+1)\int_0^\infty
t^k P\Bigl(\sup_{i\geq1} (S_i-K i)_+ >t
\Bigr)\,dt
\\
&\leq& (k+1) \sum_{i=1}^\infty\int
_0^\infty t^k P(S_i-Ki>t)
\,dt
\\
&\leq&C'\sum_{i=1}^\infty
i^{k+1/2} \exp\biggl(-i K h\biggl( -\frac
{K-1}{K}\biggr)\biggr)<
\infty
\end{eqnarray*}
for some constant $C'$ (which depends on $k$), where we have used
that\break
$h( -\frac{K-1}{K})>0$. This completes the proof.
\end{pf}

Finally, we turn to general dimension.
In our last result we combine the upper bound in Theorem~\ref{cotaenesperanza}
with a combinatorial approach to give the exact rate of convergence of the
empirical cost of optimal incomplete transportation to the uniform
distribution on the $d$-dimensional
unit cube. The result is not given in terms of expectations as in
Theorem~\ref{cotaenesperanza}. Our approach allows also to get that
type of result, but we
refrain from adding more technicalities.

\begin{Theorem} \label{tasageneral} If $P$ is the uniform distribution
on the unit
cube $[0,1]^d$, $X_1,\ldots,X_n$ are i.i.d. $P$, $P_n$ is the empirical
measure on $X_1,\ldots,X_n$ and $\alpha\in(0,1)$, then
\[
\mathcal{W}_p\bigl(\mathcal{R}_\alpha(P_n),P
\bigr)=O_P\bigl(n^{-1/d}\bigr).
\]
\end{Theorem}

\begin{pf}
For the sake of simplicity we consider the
case $d=2$.
The idea carries over smoothly to higher dimension. On the other hand,
the case $d\geq3$ follows
from known results for the usual transport; recall (\ref
{Talagrandbound}). We write $N=[\sqrt{n}]$
and $X_j=[X_{j,1},X_{j,2}]^T$, $j=1,\ldots,n$. We denote also
$B_{i}=\sharp\{j\in\{ 1,\ldots,n\}\dvtx  X_{j,1} \in(\frac{i-1} N,
\frac i N] \}$, $i=1,\ldots, N$.
The random vector $(B_{1},\ldots, B_{N})$ follows a multinomial distribution
with parameters $n$ and $(\frac1 N,\ldots, \frac1 N)$.
Given $B_i=n_i>0$,
we denote by $j^i_1,\ldots,j^i_{n_i}$ the indices $k$ such that
$X_{k,1}\in(\frac{i-1} N, \frac i N]$.
Then $X_{j_1^i,2},\ldots,X_{j_{n_i}^i,2}$
are an i.i.d. $U(0,1)$ sample. We write $P(\frac\alpha2,i)$ for the
$\frac\alpha2$-trimming of the empirical distribution on
$X_{j_1^i,2},\ldots,X_{j^i_{n_i},2}$ considered
in the proof of Theorem~\ref{cotaenesperanza}. Then we have
$E(n_i^p \mathcal{W}_p^p (P(\frac\alpha2,i), U(0,1))|B_i=n_i)\leq
C_p(\alpha/2)$. We write also $\varphi_i$ for the optimal
transportation map from $U(0,1)$ to
$P(\frac\alpha2,i)$. Then $\mathcal{W}_p^p (P(\frac\alpha2,i),
U(0,1))=\int_0^1 |x_2-\varphi_i(x_2)|^p \,dx_2$. We recall that $\varphi$
takes values on the set
$\{X_{j_1},\ldots,X_{j_{n_i}}\}$ and, with $\ell_d$ denoting
$d$-dimensional Lebesgue measure, $\ell_1 (x\dvtx  \varphi
_i(x)=X_{j_l,2})\leq\frac1 {n_i (1-\alpha/2)}$.

Next we define the map $\varphi$ on $(0,1]\times(0,1]$ as follows. If
$x=[x_1, x_2]^T$ is such that $x_1\in(\frac{i-1} N, \frac i N]$ and
$\varphi_{i}(x_2)=X_{j_l,2}$, then
$\varphi(x)=X_{j_l}$. In other words, points on the stripe $(\frac
{i-1} N, \frac i N]\times(0,1]$ are mapped to one of the
observations on that stripe, the precise one being determined by the
$\alpha/2$ trimming function on the second coordinate. Clearly, for
$x\in(\frac{i-1} N, \frac i N]$,
\[
\bigl\| x-\varphi(x)\bigr\|\leq\frac1 {N}+\bigl|x_2-\varphi_i(x_2)\bigr|.
\]
From this we get
\begin{eqnarray*}
\int_{(0,1]\times(0,1]} \bigl\| x-\varphi(x)\bigr\|^p \,dx & =&\sum
_{i=1}^N \int_{(({i-1})/N, i/N]\times(0,1]} \bigl\|
x-\varphi(x)\bigr\|^p \,dx
\\
&\leq& \frac{2^{p-1}} {N^p}+\frac{2^{p-1}} N \sum_{i=1}^N
\int_0^1 \bigl|x_2-
\varphi_i(x_2)\bigr|^p \,dx_2.
\end{eqnarray*}
Furthermore, $\ell_2 (x\dvtx  \varphi(x)=X_j)\leq\frac1 N \frac1 {n_i
(1-\frac\alpha2)}$ if $X_j\in(\frac{i-1} N, \frac i N]\times
(0,1]$. Thus, $\varphi$
maps $P$ into an $\alpha$-trimming of $P_n$ if
\[
\frac1 N \frac1 {n_i (1-\alpha/2)} \leq\frac1 {n(1-
\alpha)},\qquad i=1,\ldots, N,
\]
or, equivalently, if $\min_{1\leq i \leq N} n_i\geq\frac n N A$, with
$A=\frac{1-\alpha}{1-\alpha/2}<1$. As a consequence,
on the set $B= (\min_{1\leq i \leq N} n_i\geq\frac n N A )$,
\[
\mathcal{W}_p^p \bigl(\mathcal{R}_\alpha(P_n),P
\bigr)\leq \frac{2^{p-1}} {
N^p}+\frac{2^{p-1}} N \sum
_{i=1}^N \int_0^1
\bigl|x_2-\varphi_i(x_2)\bigr|^p
\,dx_2.
\]
Now, we note that
\begin{eqnarray*}
&&
E \Biggl(I_B \sum_{i=1}^N
\int_0^1 \bigl|x_2-\varphi_i(x_2)\bigr|^p
\,dx_2 \Big| B_1=n_1,\ldots, B_N=n_N
\Biggr)
\\
&&\qquad
\leq C_p(\alpha/ 2)I_B \sum
_{i=1}^N \frac1 {n_i^p}
\leq C(\alpha/2) \frac1 {A^p} \sum_{i=1}^N
\frac{N^p}{n^p}.
\end{eqnarray*}
But the last two displays imply that
\[
E\bigl(\mathcal{W}_p^p \bigl(\mathcal{R}_\alpha(P_n),P
\bigr)I_B\bigr)\leq\frac{2^{p-1}} {
N^p} + C_p(\alpha/2)
\frac{2^{p-1}}{A^p} \frac{N^p} {n^p}.
\]
This, toghether with the fact that $P(B^C)\to0$ as $n\to\infty$ (see
Theorem 7, page~112, in~\cite{Kolchinetal1978}) implies that
$\mathcal{W}_p^p (\mathcal{R}_\alpha(P_n),P)=O_P(n^{-p/2})$ and
completes the proof.
\end{pf}




\printaddresses

\end{document}